\documentclass[a4paper,10pt]{article}     
\usepackage{theorem}                      
\usepackage{amstex,amssymb,amscd}         

\newtheorem{theorem}{Theorem}[section]
\newtheorem{lemma}[theorem]{Lemma}
\newtheorem{proposition}[theorem]{Proposition}
\newtheorem{cor}[theorem]{Corollary}
\theorembodyfont{\rm}

\newtheorem{uremark}{Remark}

\newcommand{\CC}{{\Bbb{C}}}
\newcommand{\PP}{{\Bbb{P}}}

\newcommand{\ZZ}{{\Bbb{Z}}}

\renewcommand{\phi}{\varphi}
\newenvironment{Proof}{\begin{ProofwCaption}{Proof}}{\end{ProofwCaption}}
\newenvironment{Proof*}[1]{\begin{ProofwCaption}{{#1}}}{\end{ProofwCaption}}
\newenvironment{ProofwCaption}[1]%
  {\addvspace\theorempreskipamount \noindent{\it #1.}\rm}%
  {\qed \par \addvspace\theorempostskipamount}
\newcommand{\qedsymbol}{\mbox{$\Box$}}
\newcommand{\qed}{\quad\qedsymbol}
\begin{document}
\title{A Remark on the Geometry of Elliptic Scrolls
and Bielliptic Surfaces}
\author{C.~Ciliberto and K.~Hulek}
\date{}
\maketitle

\begin{center}
{\large \it Dedicated to the memory of F. Serrano}
\end{center}
\vspace{5mm}


\section{Introduction}
The union of two quintic elliptic scrolls in $\PP^4$
intersecting transversally along an elliptic normal quintic curve is a
singular surface $Z$ which behaves numerically like a bielliptic surface.
In the appendix to the paper \cite{ADHPR} where the equations of this
singular surface were computed, we proved that $Z$ defines a smooth point
in the appropriate Hilbert scheme and that $Z$ cannot be smoothed in
$\PP^4$. Here we consider the analogous situation in higher dimensional
projective spaces $\PP^{n-1}$, where, to our surprise, the answer depends on
the dimension $n-1$. If $n$ is odd the union of two scrolls cannot be
smoothed, whereas it can be smoothed if $n$ is even. We construct an
explicit smoothing.

\section{Elliptic scrolls}
To every elliptic normal curve $E\subset \PP^{n-1}$ of degree $n$ and
every point $P\in E$ one can associate a {\em translation scroll}
$S=S(E,P)$ by defining $S$ as the union of all secants of $E$ joining
the points $x$ and $x+P,\  (x\in E)$. If $n\ge 5$ and $P$ is not a
2-torsion point, then $S$ is a singular surface of degree $2n$ with $E$
as its singular locus. If, however, $P$ is a non-zero 2-torsion point,
then $S$ becomes a smooth scroll of degree $n$ (the secants spanned by
$x, x + P$ and $x-P, x$ coincide). Varying the point $P$ does,
of course, not define a flat
degeneration: as a scheme the general translation scrolls degenerates
to a multiplicity-$2$ scheme whose support is the smooth scroll of
degree $n$ (cf. \cite {HVdV}). In this paper we are interested in the
degree
$n$ scrolls defined by non-zero 2-torsion points. We will simply call
these scrolls {\em degree} $n$ {\em elliptic scrolls}. (If $n=5$ this
is the unique irregular smooth scroll in $\PP^4$.) Our first aim is to
determine these scrolls as abstract surfaces. Here, as in the sequel,
we shall notice a difference between the cases $n$ even and $n$ odd.\\

We shall first treat the $n$ odd case. For this purpose we fix an
elliptic curve~$F$. Recall that there is a unique $\PP^1$-- bundle over
$F$ with invariant $e=-1$ (in the sense of \cite[Chapter V]{Ha}).
It was already
observed in \cite[p. 451]{A} that this is the symmetric product
$S^2F$ where the $\PP^1-$bundle structure is given by summation
$$
\begin{array}{rll}
\pi:  &  S^2 F  &  \rightarrow F\\
& \{x,y\} &  \mapsto x+y.
\end{array}
$$
{}Fix an origin $0$ of $F$, and let $p:F\times F\rightarrow S^2F$ be the
natural projection. The curve
$$
{}F_0=p(F\times\{0\})
$$
is a section of the $\PP^1$ -- bundle $S^2F$. We choose the point $p(0,0)$
as its origin and, by abuse of notation, shall denote it again by $0$.
If $\{P_i; i=1,2,3\}$ are the non-zero 2-torsion points of $F$, then the
curves
$$
\Delta_i=\{(x,x+P_i); x\in F\} (\cong F)\quad (i=1,2,3)
$$
are mapped 2:1 under $p$ to 2-sections $F_i\subset S^2F$. As abstract
curves $F_i=F/\langle P_i\rangle$. We shall choose the point
$0_i=p(0,P_i)=p(P_i,0)$ as the origin of $F_i$. The group of 2-torsion
points of $\Delta_i$ is mapped to two points $(0_i,Q_i)$. Note that
$F_i/\langle Q_i\rangle\cong F$. Every fibre $f$ of $\pi$ intersects
$F_i$ in two points which differ by $Q_i$. We shall denote the fibre of
$\pi$ over
$P\in F$ by $f_P$ and put
$S=S^2F$. The following formulae follow immediately from the above
description:
$$
\begin{array}{ll}
(1)  &  {\cal O}_S(F_0)|_{F_0}={\cal O}_{F_0}(0)\\[1mm]
(2)  &  {\cal O}_S(F_i)={\cal O}_S (2F_0-f_{P_i})\\[1mm]
(3)  &  {\cal O}_S(F_i)|_{F_i}={\cal O}_{F_i}(0_i-Q_i)\\[1mm]
(4)  &  K_S= {\cal O}(-2F_0+f_0).
\end{array}
$$

\begin{proposition}\label{pro1.1}
Assume $n\ge 5$ odd. The line bundle ${\cal O}_S(H)={\cal
O}_S(F_0+\left(\frac{n-1}2\right)f_0)$ is very ample and embeds $S$ as
smooth surface of degree $n$ in $\PP^{n-1}$. This surface is the
translation scroll of the elliptic normal curves $F_i, i=1,2,3$ defined
by the 2-torsion points $Q_i$. Conversely every translation scroll of an
elliptic normal curve of degree $n$ by a 2-torsion point arises in this
way.
\end{proposition}

\begin{Proof}
Very ampleness of ${\cal O}_S(H)$ follows e.g.\ from
\cite[Exercise V.2.12]{Ha}.
A straightforward calculation using Riemann-Roch shows $h^0({\cal
O}_S(H))=n$. Since $H.F_i=n$ and $h^1({\cal O}_S(H-F_i))=0$ 
(the latter can be seen e.g. by Kodaira vanishing),
the
2-sections $F_i$ are mapped to elliptic normal curves of degree $n$. By
construction $S$ is then the translation scroll defined by the pair
$(F_i, Q_i)$. Conversely given any pair $(F_i, Q_i)$ consisting of an
elliptic curve and a 2-torsion point, then $F_i$ is a 2-section of
$S^2F$ with $F=F_i/\langle Q_i\rangle$ such that the rulings of $S^2F$ cut
$F_i$ in two points differing by $Q_i$.
\end{Proof}

We now turn to the case $n$ even where we assume $n\ge 6$. Let $Y\subset
\PP^{n-1}$ be a degree $n$ scroll given by a pair $(E,Q_i)$ where $E$ is
an elliptic normal curve of degree $n$ and $P_i$ a non-zero 2-torsion
point of $E$. Then $Y$ can also be constructed as follows: Embed $E$ as
a normal curve of degree $n+1$ in $\PP^n$ and let $X$ be the degree
$(n+1)$-scroll given by $E\subset \PP^n$ and the point $Q_i$. Projection
from say the origin $0\in E$ then maps $X$ to $Y$. More precisely, the
surface ${\tilde X}$ which is the blow-up of $X$ in $0$ is mapped to $Y$.
Under this projection map the fibre of $X$ over the origin $0$ is contracted,
whereas the map is bijective otherwise. In fact this is the geometric
realization of an elementary transformation of $X$.

\begin{proposition}\label{pro1.2}
The degree $n$ scroll $Y$ is smooth. It is isomorphic to the
$\PP^1$ -- bundle $\PP(\cal E)$ over the elliptic curve $F=E/\langle
Q_i\rangle$ where
${\cal E}={\cal O}_F\oplus{\cal O}_F(0-P_i)$ and $P_i$ is the image of the
2-torsion points $Q_j, j\neq i$ under the projection to $F$. The
embedding is given by the complete linear system defined by the line
bundle ${\cal O}_Y(H)={\cal O}_Y(F_0+\frac n2 f_0)$ where $F_0$ (by abuse
of notation) is the image of the section $F_0$ of $X$.
\end{proposition}

\begin{Proof}
Recall the situation on $X=S^2F$ where $E$ is the 2-section $F_i$.
There are two sections of $X$ which intersect $F_i$ (transversally) in
the point $p(0,P_i)$, namely $F_0$ and $F_{Pi}=p(F\times \{P_i\})$.
Since $F^2_0=F^2_{P_i}=F_0 F_{P_i}=1$ these sections become disjoint
sections of the $\PP^1$-bundle $\PP({\cal E})$ which is defined by the
elementary transformation with centre $p(0,P_i)$. We shall denote these
sections again by $F_0$ resp. $F_{P_i}$. Since after the elementary
transformation $F^2_0=F^2_{P_i}=0$ it follows that we can assume that
${\cal E}={\cal O}\oplus {\cal M}$ where ${\cal M}$ has degree $0$. It 
follows from the elementary transformation that the normal
bundle of $F_0$ in $\PP({\cal E})$ is isomorphic to ${\cal O}_F(0-P_i)$. This
shows the claim about ${\cal E}$. The above description of the elementary
transformation immediately gives ${\cal O}_{\PP({\cal E})}(H)={\cal
O}_{\PP({\cal E})}(F_0+\frac n2 f_0)$. Clearly $H^2=n$. Since $h^0({\cal
O}_{\PP({\cal E})}(H))=n$ and since $H$ is very ample
\cite[Exercise V.2.12]{Ha} the
proposition follows.
\end{Proof}

\begin{uremark}
Using the adjunction formula we immediately obtain the following results:

$$(5)\quad K_Y={\cal O}_Y(-2F_0)\otimes{\cal O}_Y(f_0-f_{P_i})$$

\noindent
The curve $E$ is again a 2-section of $Y$ with self-intersection number
$E^2=0$. Since $E$ and $F_0$ do not intersect we obtain

$$(6)\quad {\cal O}_Y(E)={\cal O}_Y(2F_0)$$

\noindent
and combining (5) and (6) gives

$$(7)\quad {\cal O}_Y(E)={\cal O}_Y(-K)\otimes {\cal O}_Y(f_0-f_{P_i}).$$

\noindent
Note that an analogous formula holds for $E=F_i$ on $S=S^2F$.
\end{uremark}

\section{Rigidity for $n$ odd}
We fix an elliptic normal curve $E$ in $\PP^{n-1}$ of odd degree $n$ and
two non-zero 2-torsion points $P_i\neq P_j$ on $E$. These define degree
$n$ elliptic scrolls $X_i$ and $X_j$. The union $Z=X_i\cup X_j$ of these
scrolls is a singular surface of degree $2n$ whose singular locus is the
curve $E$, which is a double curve of $Z$. Numerically $Z$ is a
bielliptic surface, its dualizing sheaf is a line bundle $\omega_Z$ with
$\omega^2_Z={\cal O}_Z$. In the case $n=5$ those surfaces were
considered in connection with abelian and bielliptic surfaces in
\cite{ADHPR}, where their equations were determined. It was
also shown \cite[appendix]{ADHPR} that they define smooth points in
their Hilbert scheme and that they are rigid, in the sense that every
small deformation of $Z$ is again of the same type. In particular, these
surfaces cannot be smoothed. In this section we shall see that this is
the same for all odd degrees, whereas, surprisingly, the situation is
very different for $n$ even.\\

Our first aim is to study the normal bundle of the degree $n$ scroll in
$\PP^{n-1}$. Let $X_i$ be one of these scrolls. Then we have the
following exact sequence for the rulings $f$ of this scroll:

$$
\begin{array}{rccccl}
0\rightarrow & N_{f/X_i}  &  \rightarrow  &  N_{f/\PP^{n-1}}  &
\rightarrow  &  N_{X_i/\PP^{n-1}}|_f\quad  \rightarrow 0\\
&  || &  &  ||  &  &\\
& {\cal O}_f  &  & (n-2){\cal O}_f(1) & &
\end{array}
$$
It follows that
$$
N_{X_i/\PP^{n-1}}|_f = (n-4){\cal O}_f(1)\oplus {\cal O}_f(2).
$$
The degree 2 subbundle is uniquely determined and varying $f$ we obtain
a line subbundle
$$
{\cal L} = \pi^*(\pi_*N_{X_i/\PP^{n-1}}(-2))(2)\subset N_{X_i/\PP^{n-1}}
$$
and as in \cite{HVdV} one proves that ${\cal L}=K^{-1}_{X_i}$.
Thus we have an exact sequence

$$(8)\quad 0\rightarrow K^{-1}_{X_i} \rightarrow N_{X_i/\PP^{n-1}}\rightarrow
Q\rightarrow 0$$

\noindent where $Q$ is a vector bundle of rank $(n-4)$ with
$Q|_f=(n-4){\cal O}_f(1)$.

\begin{lemma}\label{lem2.1}
\begin{enumerate}
\item[${\rm(i)}$]
$h^0(N_{X_i/\PP^{n-1}})=n^2$
\item[${\rm(ii)}$]
$h^j(N_{X_i/\PP^{n-1}})=0 \mbox{ for } j\ge1.$
\end{enumerate} 
\end{lemma}

\begin{Proof}
Since $h^j(K^{-1}_{X_i})=0$ for all $j$, it follows from sequence (8)
that the claim is equivalent to $h^0(Q)=n^2$ and $h^j(Q)=0$ for $j\ge
1$. The defining sequences for $Q$ and the normal bundle
$N_{X_i/\PP^{n-1}}$ together with Riemann-Roch give
$$
\chi(Q)=\frac 12 (c^2_1(Q)-2c_2(Q))+\frac 12
c_1(Q)(-K_{X_i})+(n-4)\chi({\cal O}_X)=n^2
$$
Hence it is enough to prove that $h^j(Q)=0$ for $j\ge 1$. By Serre
duality $h^2(Q)=h^0(Q^{\vee} \otimes K_{X_i})=0$ since $Q^{\vee}\otimes
K_{X_i|f}=(n-4){\cal O}_f(-3)$. To prove vanishing of $h^1(Q)$ we first
remark that $Q(-1)$ is trivial on the fibres $f$ and hence
$Q(-1)=\pi^*{\cal F}$ where ${\cal F}$ is a rank $n-4$ bundle on $F$.
Since $T_{\PP^{n-1}}(-1)$ is generated by global sections the same is
true for $N_{X_i/\PP^{n-1}}(-1)$ and hence also for $Q(-1)$ and ${\cal
F}$. But now, using the classification of vector bundles on elliptic
curves \cite{A} it follows that $h^1({\cal F}(D))=0$ for every divisor
$D$ on $F$ of positive degree. Recall that $H=F_0+\left(
\frac{n-1}2\right) f_0$. Let $Q'=Q(-1)\otimes \left(
\frac{n-1}2\right) f_0$. Then $h^1(Q')=h^1({\cal F}(\frac{n-1}2
0))=0$. Finally $h^1(Q)=0$ follows from the exact sequence
$$
0\rightarrow Q'\rightarrow Q\rightarrow Q|_{F_0}\rightarrow 0
$$
since $Q|_{F_0}={\cal F}(\frac{n+1}2 f_0)$ and $h^1({\cal F}(\frac{n+1}2
f_0))=0.$
\end{Proof}

\begin{lemma}\label{lem2.2}
$h^j(N_{X_i/\PP^{n-1}}(-F_i))=0$ for all $j$.
\end{lemma}

\begin{Proof}
Twisting (8) with ${\cal O}_{X_i}(-F_i)$ we obtain the exact sequence
$$
0\rightarrow K^{-1}_{X_i}(-F_i)\rightarrow
N_{X_i/\PP^{n-1}}(-F_i)\rightarrow Q(-F_i)\rightarrow 0.
$$
The line bundle $K^{-1}_{X_i}(-F_i)={\cal O}_{X_i}(f_{P_i}-f_0)$ has no
cohomology. Since for the restriction to a ruling
$Q(-F_i)|_f=(n-4){\cal O}_f(-1)$ it follows that
$h^0(Q(-F_i))=h^2(Q(-F_i))=0$. Finally we obtain
$-h^1(Q(-F_i))=\chi(Q(-F_i))=0.$
\end{Proof}

We now turn to the normal bundle $N_{Z/\PP^{n-1}}$ of $Z$.

\begin{proposition}\label{pro2.3}
\begin{enumerate}
\item[${\rm(i)}$]
$h^0(N_{Z/\PP^{n-1}})=n^2$
\item[${\rm(ii)}$]
$h^j(N_{Z/\PP^{n-1}})=0$ for $j\ge 1$.
\end {enumerate}
\end{proposition}

\begin{Proof}
The line bundle
$$
T=N_{E/X_i}\otimes N_{E/X_j}={\cal O}_E(20-Q_i-Q_j)
$$
is a non trivial 2-torsion bundle. As in \cite{CLM} we have the following
exact sequences

$$(9)\quad 0\rightarrow N_{X_i/\PP^{n-1}}\rightarrow
N_{Z/\PP^{n-1}}|_{X_i}\rightarrow T\rightarrow 0$$

$$(10)\quad 0\rightarrow N_{X_i/\PP^{n-1}}(-E)\rightarrow
N_{Z/\PP^{n-1}}|_{X_i}\otimes {\cal O}_{X_i}(-E)\rightarrow T\otimes {\cal
O}_{X_i}(-E)\rightarrow 0$$

$$(11)\quad 0\rightarrow N_{Z/\PP^{ n-1}}|_{X_j}\otimes {\cal
O}_{X_j}(-E)\rightarrow N_{Z/\PP^{n-1}}\rightarrow
N_{Z/\PP^{n-1}}|_{X_i}\rightarrow 0.$$

By formula (3) $T\otimes{\cal O}_{X_i}(-E)={\cal O}_E(0-Q_j)$. Together
with Lemma \ref{lem2.2} it follows from (10) that
$$
h^j (N_{Z/\PP^{n-1}}|_{X_i}\otimes {\cal O}_{X_i}(-E))=0 \mbox{ for all
} j.
$$
{}From sequence (9) and Lemma \ref{lem2.1} we obtain
$$
h^0(N_{Z/\PP^{n-1}}|_{X_i})=n^2,\quad h^j(N_{Z/\PP^{n-1}}|_{X_i})=0
\mbox { for } j\ge 1.
$$
The result now follows from sequence (11).
\end{Proof}

\begin{theorem}\label{theo2.4}
(Rigidity) The component of the Hilbert scheme of surfaces containing $Z$
is smooth of dimension $n^2$ at [Z]. Every small deformation of $Z$
is again a union of two degree $n$ elliptic scrolls intersecting
transversally along an elliptic normal curve.
\end{theorem}

\begin{Proof}
The statement about smoothness and the dimension of the Hilbert scheme
follows immediately from Proposition \ref{pro2.3}.
Let $X(2,n)$ be the modular curve parametrizing elliptic curves with a
level $n$ structure and a non-zero 2-torsion point. Every point of
$X^0(2,n)=X(2,n)-\mbox{\{cusps\}}$ gives rise to a union $Z$ of
two degree
$n$ elliptic scrolls. Indeed the elliptic curves with level $n$ structure
are in 1:1 correspondence with Heisenberg invariant elliptic normal
curves in $\PP^{n-1}$. Given a non-zero 2 torsion point we have exactly
two other such points. We can use these two points to construct~$Z$.
Conversely every point $Z$ arises in this way up to a change of
coordinates. Let ${\cal H}$ be the component of the Hilbert scheme
containing a given surface $Z$. Then we have a natural map
$$
\Phi: X^0(2,n)\times \operatorname{PGL } (n,\CC)\rightarrow {\cal H}.
$$
Since every elliptic normal curve has a finite automorphism group this
map is finite and hence surjective in a neighbourhood of [Z].
\end{Proof}

\begin{uremark}
Of course there exist global deformations of $Z$ which are not a union
of two degree $n$ elliptic scrolls. E.g. $Z$ can degenerate into
non-reduced union of $n$-planes. For a discussion of possible
degenerations in the case of $\PP^4$, i.e. $n=5$ see
\cite[section 9]{ADHPR}.
\end{uremark}

\section{Smoothing for $n$ even}

The case $n$ even is subtly different from the case $n$ odd, as can
already be seen in the computation of the cohomology of the normal
bundle of the union of two scrolls. Again we fix an elliptic normal
curve $E$ in $\PP^{n-1}$ of degree $n$ and two non-zero 2-torsion points
$P_i\neq P_j$ defining degree $n$ elliptic scrolls $X_i$ and $X_j$. Let
$Z=X_i\cup X_j$. As before we find for the normal bundle of $X_i$ an
exact sequence
$$
0\rightarrow K^{-1}_{X_i}\rightarrow N_{X_i/\PP^{n-1}}\rightarrow
Q\rightarrow 0
$$
with $Q|_f=(n-4){\cal O}_f(1)$. In this case, however,
$h^0(K_{X_i}^{-1})=h^1(K_{X_i}^{-1})=1$. Nevertheless the arguments of 
Lemma \ref{lem2.2} still go through and give

$$(12)\quad  h^j(N_{X_i/\PP^{n-1}}(-E))=0 \mbox{ for all } j.$$

We also have an exact sequence

$$(13)\quad 0\rightarrow N_{X_i/\PP^{n-1}}(-E)\rightarrow
N_{X_i/\PP^{n-1}}\rightarrow N_{X_i/\PP^{n-1}}|_E\rightarrow 0.$$

Since $N_{X_i/\PP^{n-1}}(-1)$ is globally generated we can conclude as
in the proof of Lemma \ref{lem2.1} that $h^j(N_{X_i/\PP^{n-1}}|_E)=0$
for $j\ge 1$. But then sequence (13) together with Riemann-Roch
(numerically the cases $n$ odd and $n$ even behave in exactly the same
way) gives

$$
(14)\quad  h^0( N_{X_i/\PP^{n-1}})=n^2,\quad h^j
(N_{X_i/\PP^{n-1}})=0\mbox { for } j\ge 1
$$
which is as in the degree $n$ odd case. The main difference between the
two cases lies in the fact that $N_{E/X_i}=N_{E/X_j}={\cal O}_E$ and
hence

$$(15)\quad  T=N_{E/X_i} \otimes N_{E/X_j}={\cal O}_E.$$

It now follows from sequences (9) and (11) together with formula (14) that
$h^1(N_{Z/\PP^{n-1}}|_{X_i})=1$. Since it follows by sequence (10) and by
(12) that $h^2(N_{Z/\PP^{n-1}}|_{X_j}\otimes {\cal O}_{X_j}(-E))=0$ we 
find that $h^1(N_{Z/\PP^{n-1}})>0$, contrary to the degree $n$ odd case. 
Moreover
sequences (9)--(11) show that $h^0(N_{Z/\PP^{n-1}})=n^2+1$ or $n^2+2$ and
$h^1(N_{Z/\PP ^{n-1}})=1$ or $2$ respectively. In fact we shall see later 
(Corollary \ref {cor3.3}) that $h^0(N_{Z/\PP^{n-1}})=n^2+1$.\\

We now want to construct an explicit embedded smoothing of the singular
surface $Z$ to a bielliptic surface. Since $\omega^2_Z={\cal O}_Z$ it is
natural to look at bielliptic surfaces of type 1) or 2) in the Bagnera-de
{}Franchis list \cite [List VI.20]{B}. It is easy to see by Reider's method 
(cf.\cite {Se}) that bielliptic surfaces of type 1) cannot be embedded in 
$\PP^{n-1}$ for $n\le 8$. Hence we shall now turn our attention to 
bielliptic surfaces of type 2). Recall that these
surfaces are of the form $S=E\times F/G$ where $G=\ZZ_2\times \ZZ_2$ acts 
on $E$ by translation with 2-torsion points and on $F$ by $x\mapsto -x,
x\mapsto x+\varepsilon, \varepsilon$ a 2-torsion point of $F$. We shall first
show that these surfaces can be embedded as surfaces of degree $2n$ in
$\PP^{n-1}$. This will then give us the right idea for the construction
of the degenerations.

\begin{proposition}\label{pro3.1}
Every bielliptic surface $S$ of type 2) can be embedded as a linearly
normal surface of degree $2n$ in $\PP^{n-1}(n\ge 6)$.
\end{proposition}

\begin{Proof}
Let $\pi:E\times F\rightarrow S$ be the projection map and set
$A=\pi(E), B=\pi(F)$. Then $A.B=4$.
By \cite[Proposition 1.7]{Se} the element $B/2$ is in
$\operatorname{NS}(S)$ and we can consider the divisor
$$
H=A+\frac n 4 B.
$$
(Since $n$ is even this is indeed a divisor on $S$). Then $H^2=2n,
H.A=n$ and $H.B=4$. It is easy to check that $H$ is ample and
Riemann-Roch together with Kodaira vanishing gives $h^0({\cal
O}_S(H))=n$. It is a straightforward application of Reider's theorem to
prove that $H$ is very ample. Hence the complete linear system defined
by $H$ embeds $S$ as a linearly normal surface of degree $2n$ in
$\PP^{n-1}$.
\end{Proof}

\begin{uremark}
The line bundle $\pi^*{\cal O}_S(H)$ has degree $n$ on $E$ and degree
$4$ on $F$.
\end{uremark}

\begin{theorem}\label{theo3.2}
(Smoothing) Let $Z=X_i\cup X_j$ be a union of two degree $n$ elliptic
scrolls in $\PP^{n-1}(n\ge 6$, even). Then there exists a flat family of
surfaces $(Z_t)_{t\in T}$ in
$\PP^{n-1}$ such that $Z_0=Z$ and $Z_t$ for $t\neq 0$ is a
linearly normal smooth bielliptic surface of degree $2n$.
\end{theorem}

\begin{Proof}
We fix the elliptic curve $E=X_i\cap X_j$ and the two non-zero 2-torsion
points $P_i$ and $P_j$ which define $X_i$ and $X_j$. Let
$F=\CC/(\ZZ+\ZZ\tau)$ be another elliptic curve which we consider
variable. Let $S(4)\rightarrow X(4)$ be the Shioda modular surface of
level 4. We consider the family ${\cal F}=(F_t)_{t\in T}$ where $t$ varies
in some neighbourhood of a cusp of $X(4)$, say $i\infty$,
where $t=e^{2\pi i {\tau}/4}$. Then $F_0$ is a
4-gon of rational curves and the 2-torsion points $Q_0=0, Q_1=1/2,
Q_2=\tau/2$ and $Q_3=(1+\tau)/2$ of $F_t=\CC/(\ZZ+\ZZ\tau)$ define 4
sections of ${\cal F}$ which intersect the singular fibre $F_0$ as
indicated below

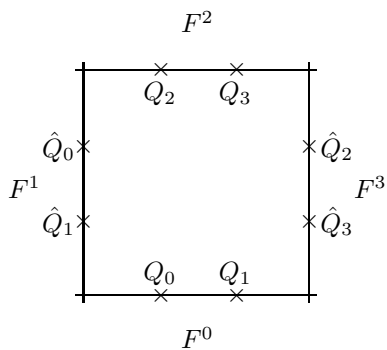
\begin{figure}[h]
$$
\unitlength1cm
\begin{picture}(5,4.5)
\put(0.9,0.7){\line(1,0){3.2}} \put(0.9,3.7){\line(1,0){3.2}}
\put(1.0,0.6){\line(0,1){3.2}} \put(4.0,0.6){\line(0,1){3.2}}
\put(0.865,1.6){$\times$} \put(0.865,2.6){$\times$}
\put(3.865,1.6){$\times$} \put(3.865,2.6){$\times$}
\put(1.9,3.615){$\times$} \put(2.9,3.615){$\times$}
\put(1.9,0.615){$\times$} \put(2.9,0.615){$\times$}
\put(1.8,0.9){$Q_0$} \put(2.8,0.9){$Q_1$}
\put(1.8,3.3){$Q_2$} \put(2.8,3.3){$Q_3$}
\put(0.45,1.55){${\hat Q}_1$} \put(0.45,2.55){${\hat Q}_0$}
\put(4.15,1.55){${\hat Q}_3$} \put(4.15,2.55){${\hat Q}_2$}
\put(0,2){$F^1$} \put(4.6,2){$F^3$}
\put(2.3,0){$F^0$} \put(2.3,4.2){$F^2$}
\end{picture}
$$
\caption{Singular fibre of $\cal F$}
\end{figure}

\noindent The action of the 2-torsion points on smooth fibres extends to
an action on~${\cal F}$. Similarly the involution
$ \iota : x\mapsto-x$
on smooth fibres extends to an involution on ${\cal F}$. We
denote the section of ${\cal F}$ given by $Q_2$ by $\varepsilon$. Then
$\varepsilon$ acts on $F_0$ by rotation with $180^{\circ}$, i.e. it
identifies
$F^0$ and $F^2$, resp. $F^1$ and $F^3$. The involution $\iota$ interchanges
$F^1$ and $F^3$, resp.\ induces involutions on $F^0$ and $F^2$ with fixed
points $Q_0, Q_1$ and $Q_2, Q_3$.\\

We now consider the product ${\cal X}=E\times {\cal F}$ which is naturally
fibred over $T$ with fibre $X_t=E\times F_t$. We define an action of
$G=\ZZ_2\times \ZZ_2$ on ${\cal X}$ as follows: The element
$g_1=(1,0)$ acts on $E$ by $x\mapsto x + P_i$ and on $F$ by $x\mapsto
x+\varepsilon$, whereas $g_2=(0,1)$ acts on $E$ by $x\mapsto x+P_j$ and on
$F$ by $x\mapsto -x$. Then $g_1 g_2=(1,1)$ acts on $E$ by $x\mapsto
x+(P_i+P_j)$ and on $F$ by $x\mapsto-x+\varepsilon$. The total space ${\cal
F}$ and hence ${\cal X}$ is smooth and $G$ acts freely on ${\cal
X}$. Then ${\cal Z}={\cal X}/G$ is smooth and $Z_t=X_t/G$ is
a bielliptic surface of type 2) for $t\neq 0$. The singular surface $Z_0$
has the following properties:

\begin{itemize}
\item
$Z_0$ consists of two components $Z_0^0$ and $Z_0^1$, namely the images
of $E\times F^0$ (resp.\ $E\times F^2$) and $E\times F^1$ (resp.\
$E\times F^3$).

\item
The singular locus $E\times \operatorname{ Sing } F_0$ of $X_0$ is mapped
to an irreducible curve isomorphic to $E$ (and again denoted by $E$).
This curve $E=Z^0_0\cap Z^1_0$ is the singular locus of $Z_0$.

\item
$Z^0_0$ and $Z^1_0$ are $\PP^1$ -- bundles over the elliptic curves
$E/\langle P_i\rangle$, resp. $E/\langle P_j\rangle$ and the singular
curve $E$ is a bisection of both $\PP^1$ -- bundles.

\item
The curves $E\times \{Q_i\}$, $i=0,\ldots , 3$ are mapped to two disjoint
sections $C^0_0$ and $C^0_1$ of $Z^0_0$ with $(C^0_0)^2=(C^0_1)^2=0$.

\item
Similarly we can consider the sections of ${\cal F}$ given by the points
${\hat Q}_i=Q_i+ {\tau}/4$. These sections again intersect $F_0$ in 4
points, which this time lie on $F^1$ and $F^3$ (see again figure 1).
These curves $E\times \{{\hat Q}_i\}$ map to two sections $C^1_0$ and
$C^1_1$ on $Z^1_0$ with $(C^1_0)^2=(C^1_1)^2=0$.
\end{itemize}

The next step is to construct a suitable line bundle ${\cal L}$ on ${\cal
X}$ which descends to ${\cal Z}$. First consider the degree $n$ line
bundle ${\cal L}_0={\cal O}_E(n0)$ on $E$. Then the group $G$, which
operates on $E$ by translation with 2-torsion points leaves ${\cal L}_0$
fixed as a line bundle. However, if we want to lift the action of $G$ to
the bundle ${\cal L}$  itself we might have to extend the group $G$
depending on the commutator of lifts $g^{{\cal L}_0}_1$ and $g^{{\cal
L}_0}_ 2$ of $g_1$ and $g_2$ to ${\cal L}_0$. By general Heisenberg theory
$$
\left [g^{{\cal L}_0}_1, g^{{\cal L}_0}_2\right ]=\left( e^{2\pi
i/n}\right)^{n^2/4}=e^{2\pi i n/4}
$$
which is either $1$ or $-1$ depending on whether $n\equiv 0\mod 4$ or not.
Hence if ${n\equiv 0\mod 4}$ then the action of $G$ on $E$ lifts to an
action
on ${\cal L}_0$,whereas if $n\equiv 2\mod 4$ we have to extend $G$
to the level
2 Heisenberg group $H$ which is a central extension
$$
1\rightarrow \{\pm 1\}\rightarrow H\rightarrow G\rightarrow \{0\}.
$$
Next we consider the sections $D_i$, resp.\ $\hat {D}_i$ of ${\cal F}$
given by the points $Q_i$, resp. ${\hat Q}_i$. Let ${\cal L}_1={\cal
O}_{\cal F}(D_0+D_2+{\hat D}_0+{\hat D}_2)$. We claim that the
action of $G$ on ${\cal F}$ lifts to~${\cal L}_1$, i.e. that $\left [
g_1^{{\cal L}_1}, g_2^{{\cal L}_2}\right ]=1$. It is enough to check this
on a general fibre $F_t$ of~${\cal F}$. For this let $s$ be a section of
${\cal L}_1$ vanishing on $D_0+D_2+{\hat D}_0+{\hat D}_2$. Since this
divisor is invariant under $G$ it follows that both $g_1^{{\cal L}_1}$
and $g_2^{{\cal L}_2}$ map $s_t=s|_{F_t}$ to a multiple of itself and
hence commute. If $n\equiv 0\mod 4$ we can set ${\cal L}={\cal
L}_0\boxtimes {\cal L}_1$. Then $G$ acts on ${\cal L}$ and since $G$
acts freely on ${\cal X}=E\times {\cal F}$ we obtain a line bundle
${\bar {\cal L}}={\cal L}/G$ on ${\cal Z}$.
In the case $n\equiv 2\mod 4$ we
must replace ${\cal L}$ by some suitable other line bundle. Let ${\cal
M}_1={\cal O}_{\cal F}(D_0-D_1)$. Then ${\cal M}_1$ is invariant under
$G$, but $G$ does not lift to an action on ${\cal M}_1$. In fact we claim
that $\left [ g_1^{{\cal M}_1}, g_2^{{\cal M}_2}\right]=-1$. For this
consider the function
$$
f(\tau,z)=\frac{\vartheta(\tau,z+\frac 12 (\tau+1))}{\vartheta(\tau,
z+\frac 12\tau)}
$$
where
$\vartheta(\tau, z) = \sum\limits_{n\in\ZZ} e^{2\pi
i(\frac 12 n^2\tau+nz)}$
is the standard theta function.
This is a meromorphic section of ${\cal M}_1$ at least for $t=e^{2\pi
i/4}\neq 0$. The claim then
follows from the identity
$$
\frac{\vartheta(\tau,-z-\frac
12(\tau+1)+\frac{\tau}2)}{\vartheta(\tau,-z-\frac{\tau}2+\frac{\tau}2)}=-
\frac{\vartheta(\tau,-z-\frac
12(\tau+1)-\frac{\tau}2)}{\vartheta(\tau,-z-\frac{\tau}2-\frac{\tau}2)}.
$$
which follows immediately from \cite [pp. 49,50]{I}, formulae
($\Theta$1)-($\Theta$3). Now consider ${\cal L}={\cal L}_0 \boxtimes
({\cal L}_1\otimes{\cal M}_1)$. This time the action of $G$ on ${\cal X}$
lifts a priori to an action of
$H$ on ${\cal L}$. But by construction the centre of $H$ acts by
$(-1)^2=1$, i.e. trivially. Hence we obtain again an action of $G$ on
${\cal L}$ and we can take ${\bar{\cal L}}={\cal L}/G$.\\

It remains to verify that  ${\bar{\cal L}}$ has the desired
properties. Let ${\bar{\cal L}}_t= {\bar{\cal L}}|_{Z_t}$. For $t\neq
0$ by Proposition \ref {pro3.1}
${\bar{\cal L}}_t$ embeds $Z_t$ as a linearly normal
bielliptic surface (which by construction is of type 2)\ ). We have to
show that  ${\bar{\cal L}}_0$ embeds $Z_0$ as the union of the two
scrolls $X_i$ and $X_j$. Let  ${\bar{\cal L}}_0^i= {\bar{\cal
L}}_0|_{Z^i_0}$ for $i=0,1$. By construction  ${\bar{\cal L}}_0^i$ has
degree~$n/2$ on the sections $C_0^i$ and $C^i_1$, degree $1$ on the
rulings and degree $n$ on the bisection $E$. Hence  ${\bar{\cal
L}}_0^i\equiv {\cal O}_{Z_0^i}(C^i_0+\frac n2 f)$. Thus $h^0( {\bar{\cal
L}}_0^i)=n$ and the restriction map $H^0(Z^i_0, {\bar{\cal
L}}_0^i)\rightarrow H^0(E, {\bar{\cal L}}_0^i|_E)$ is an isomorphism. In
particular we find that
$$
h^0(Z_0,{\bar{\cal L}}_0)=h^0(Z^0_0,{\bar{\cal L}}^0_0)+h^0(Z^1_0,
{\bar{\cal L}}^1_0)-h^0(E,{\bar{\cal L}}_0|_E)=n=h^0(Z_t,{\bar{\cal
L}}_t).
$$
Moreover the restriction map $H^0(Z_0,{\bar{\cal L}}_0)\rightarrow
H^0(Z^i_0,{\bar{\cal L}}^i_0)$ is an isomorphism and  ${\bar{\cal L}}_0$
embeds each of the component $\bar{Z}^i_0$ as a degree $n$ elliptic
scroll. By construction the image scrolls are the translation scrolls of
the embedded elliptic normal curve $E$ defined by the 2-torsion
points $P_i$ and $P_j$. This gives the claim.
\end{Proof}

\begin{uremark}
The difference between the cases $n\equiv 0\mod 4$ and $n\equiv 0\mod 2$
is easily understood in terms of the geometry of bielliptic surfaces. In
the first case $\frac n4B$ is an integer multiple of $B$ and hence
effective. In the second case $B/2$ is an element of the Neron-Severi
group of $S$, but is not effective.

\end{uremark}

\begin{cor}\label{cor3.3}
If $n$ is even, then
\begin{enumerate}
\item[{\rm(i)}]
$h^0(N_{Z/\PP^{n-1}})=n^2+1$ (and hence $h^1(N_{Z/\PP^{n-1}})=1 $),
\item[{\rm(ii)}]
the Hilbert scheme containing the surface $Z$ is smooth at $[Z]$ where it
has dimension $n^2+1$.
\end{enumerate}
\end{cor}

\begin{Proof}
The bielliptic surfaces of type 2) define a component of the Hilbert scheme
containing $Z$ which is of dimension at least $n^2+1$. Hence, if we can prove
(i), then assertion (ii) is an immediate consequence. In view of our earlier
computations it is, therefore, enough to exclude the case
$h^0(N_{Z/\PP^{n-1}})=n^2+2$. Consider the diagram

$$
\begin{array}{rc}
& 0{\ \ }\\
& \downarrow{\ }\\
& H^0(N_{X_i/\PP^{n-1}})\\[1mm]
& \downarrow^{\beta}\\[1mm]
0\rightarrow H^0(N_{Z/\PP^{n-1}}|_{X_j}\otimes {\cal O}_{X_j}(-E))\rightarrow
H^0(N_{Z/\PP^{n-1}})
\stackrel{\alpha}{\rightarrow} & H^0(N_{Z/\PP^{n-1}}|_{X_i})\\[1mm]
& \downarrow{\ }\\[1mm]
& H^0(T).{\ }
\end{array}
$$

It follows from (12) and sequence (10) that $h^0(N_{Z/\PP^{n-1}}|_{X_j}
\otimes{\cal O}_{X_j}(-E))=1$. Hence, if $h^0(N_{Z/\PP^{n-1}})=n^2+2$ then, 
since $h^0(N_{Z/\PP^{n-1}}|_{X_j})=n^2+1$ (by (14) and sequence (9)), the map 
$\alpha$ must be surjective. In particular im $(\alpha) \supset \mbox{ im } 
(\beta)$. On the other hand the sequence

$$
0\rightarrow \begin{array}[t]{c}N_{E/X_i}\\\Vert\\{\cal O}_E\end{array} 
\rightarrow N_{E/\PP^{n-1}} \rightarrow N_{X_i/\PP^{n-1}}|_E\rightarrow 0
$$
gives rise to a diagram

$$
\begin{array}{cccccc}
H^0(N_{E/X_i}) &  \rightarrow  H^0(N_{E/\PP^{n-1}}) 
\rightarrow  &H^0(N_{X_i/\PP^{n-1}}|_E)&   \rightarrow  &H^1(N_{E/X_i})&
\rightarrow  0.\\[1mm]
  ||&  &\uparrow{\cong}& & || &\\[1mm]
 \CC&&H^0(N_{X_i/\PP^{n-1}}).&&\CC
\end{array}
$$
Here the vertical isomorphism is a consequence of (12). In particular we can
find a section $s\in H^0(N_{X_i/\PP^{n-1}})$  which does not lift to
$H^0(N_{E/\PP^{n-1}})$. We claim that such a section cannot be in the image of
$\alpha$. To see this, assume that there exists a section ${\tilde s}\in
H^0(N_{Z/\PP^{n-1}})$ with $\alpha({\tilde s})=\beta(s)$. Then $s$ and
${\tilde s}$ define first order deformations ${\cal X}_i$ and ${\cal Z}$ of 
$X_i$ and
$Z$ over $\mbox{ Spec } (\CC[\epsilon]/\epsilon^2)$ such that ${\cal X}_i
\subset {\cal Z}$. A straightforward local calculation then shows that 
${\cal Z}={\cal X}_i \cup {\cal X}_j$ where ${\cal X}_j$ is a first order 
deformation of $X_j$. Moreover ${\cal X}_i \cap {\cal X}_j={\cal E}$ is a 
first order deformation of $E$. In particular ${\cal X}_i$ contains a first 
order deformation of $E$ which contradicts our choice of the section s. This 
proves the claim. 
\end{Proof}

 \bigskip
 \parindent=0pt
 \medskip

Ciro Ciliberto,\par
Dipartimento di Matematica, Universita di Tor Vergata (Roma II),\par
Via Fontanile di Carcaricola, I--00133 Roma (Italy)\par
{\em E-mail address:} cilibert\symbol{64}axp.mat.utovrm.it
\vspace{5mm}

Klaus Hulek,\par
Institut f\"ur Mathematik, Universit\"at Hannover,\par
Postfach 6009, D--30060 Hannover (Germany)\par
{\em E-mail address:} hulek\symbol{64}math.uni-hannover.de\par

\end{document}